\documentclass[12pt]{article}

\usepackage{amssymb}
\usepackage{amsmath}
\usepackage{amscd}
\usepackage{setspace}
\numberwithin{equation}{section}

\newtheorem{theorem}{Theorem}[section]
\newtheorem{corollary}{Corollary}[section]

\begin{document}

\title{An Explicit Solution to the Chessboard Pebbling Problem}

\author{
Qiang Zhen\thanks{
Corresponding author. Department of Mathematics and Statistics,
University of North Florida, 1 UNF Dr, Bldg 14/2731, Jacksonville, FL 32224-7699, USA. 
{\em Email:} q.zhen@unf.edu. Tel: (904) 620-2042. Fax: (904) 620-2818.}
\and
and
\and
Charles Knessl\thanks{
Department of Mathematics, Statistics, and Computer Science,
University of Illinois at Chicago, 851 South Morgan Street (M/C 249),
Chicago, IL 60607-7045, USA.
{\em Email:} knessl@uic.edu. This author's work was partly supported by NSA grant H 98230-08-1-0102.
}}
\date{September 24, 2010 }
\maketitle

\begin{abstract}
\noindent  We consider the chessboard pebbling problem analyzed by Chung, Graham, Morrison and Odlyzko \cite{Chu}. We study the number of reachable configurations $G(k)$ and a related double sequence $G(k,m)$.   Exact expressions for these are derived, and we then consider various asymptotic limits.

\vspace{5mm}
\end{abstract}

\section{Introduction}
The following problem has attracted some attention recently. We begin with an infinite chessboard, which consists of the lattice points $\{(i, j):\, i, j \ge 0\}$ in the first quadrant. We refer to an individual lattice point as a ``cell". We start with a single pebble placed at $(0, 0)$. The first ``step" consists of removing this pebble and placing two pebbles at the cells $(0, 1)$ and $(1, 0)$. At each subsequent step we remove a pebble at cell $(i, j)$ and place two pebbles at cells $(i+1, j)$ and $(i, j+1)$, provided that the latter two cells are unoccupied. We consider all possible choices of $(i,j)$. After $k$ steps, there will be a total of $k + 1$ pebbles on the board, in various arrangements or configurations. 

We let $G(k)$ be the total number of reachable configurations with $k$ pebbles. We define the level sets $L(l)=\{(i,j) : i+j=l\}$, so that $\cup_{l\ge 0}L(l)$ is the entire quadrant. The original problem, posed by Kontsevich \cite{Kon}, was to show that $L(1)\cup L(2)\cup L(3)$ is unavoidable, in that such a set must contain at least one pebble for any reachable configuration. A partial analysis of this fact was published thereafter by Khodulev \cite{Kho}. The first complete proof was given by Chung, Graham, Morrison and Odlyzko \cite{Chu}. It was also shown in \cite{Chu} that $L(1)\cup L(2)$  is an unavoidable set, as well as certain properties relating to the number of unavoidable sets with $k$ pebbles, including the geometric growth rate of this quantity as $k\to\infty$.  In \cite{Kne}, Knessl obtained some further asymptotic properties relating to the enumeration of unavoidable sets. Various extensions of this problem are studied by Eriksson \cite{Eri} and Warren \cite{War}.

Suppose we allow more than one pebble per cell and start with an initial configuration of one pebble in cells $(0,m+1)$ and $(m+1, 0)$ and two pebbles in each of the cells $(1,m), (2,m-1), . . . , (m-1, 2), (m, 1)$. Thus there are a total of $2m+2$ pebbles in the level set $L(m+1)$, and we assume that $L(M)$ are empty for $M > m +1$. Again the pebble at $(i, j)$ can only be moved if the cells $(i+1, j)$ and $(i, j+1)$ are empty. Let the number of reachable configurations corresponding to this starting arrangement be denoted by $G(k,m)$. In \cite{Chu}, a recurrence relation for $G(k,m)$ is derived, and it is shown that, for $k\ge 2$, $G(k, 0) = G(k)$. It is also established that as $k\to\infty$, $G(k)\sim c_*a^k$ where $a = 2.321642199494\cdots$ and $c_*=0.12268707\cdots$. These constants are characterized in terms of a continued fraction representation of the generating function of $G(k)$. In \cite{Knessl}, Knessl obtained a more explicit analytic characterization of the growth rate $a$ of the number of reachable configurations, and showed that this constant can be obtained by solving a transcendental equation that involves two series, each resembling Jacobi elliptic functions. Other asymptotic properties of $G(k,m)$ for $k$ and/or $m \to\infty $ are also established in \cite{Knessl}.

In this paper, we derive an exact expression for the number of reachable configurations $G(k)$ and $G(k,m)$.  Using the exact expression, we recover the asymptotic results given in \cite{Knessl}, and obtain more explicit analytic expressions for each asymptotic scale. These asymptotic scales are $k\to\infty$ with $m=O(1)$; $k$ and $m\to\infty$ simultaneously with $2k/m^2>1$; and $k,m\to\infty$ with $l=k-m(m+5)/2=O(1)$. 

The paper is organized as follows. In Section 2, we state the basic equations and summarize the main results for $G(k,m)$ in Theorem 2.1 and 2.2, and for $G(k)$ in Corollary 2.1 and 2.2. In Section 3, we provide brief derivations.

\section{Summary of results}
It is shown in \cite{Chu} that $G(k,m)$ satisfies the following recurrence equations:
\begin{eqnarray}
G(k,0)&=&2G(k-1,0)+G(k,1)+\delta(k,2)\label{s2_G0}\\
G(k,1)&=&G(k-3,0)+2G(k-2,1)+G(k-1,2)+G(k-4,1)\label{s2_G1}\\
G(k,m)&=&G(k-m-2,m-1)+2G(k-m-1,m)+G(k-m,m+1),\; m\ge 2.\label{s2_Gk}
\end{eqnarray}
Here $\delta(i,j)$ is the Kronecker delta symbol.  As shown in \cite{Knessl}, the boundary condition (\ref{s2_G0}) can be replaced by
\begin{equation}\label{s2_G0ic}
G(k,0)=2^{k-2}+\sum_{l=1}^{k}2^{k-l}G(l,1),\; k\ge 2.
\end{equation}
Using (\ref{s2_G0ic}) in (\ref{s2_G1}), $G(k,0)$ can be eliminated, which leads to
\begin{equation}\label{s2_G1ic}
G(k,1)=2G(k-2,1)+G(k-1,2)+G(k-4,1)+2^{k-5}+\sum_{l=1}^{k-3}2^{k-l-3}G(l,1),\;k\ge 5.
\end{equation}
We introduce a generation function of $G(k,m)$, with
\begin{equation}\label{s2_Vm}
V_m(z)=(-1)^mz^{-m}\sum_{k=0}^\infty G(k,m)z^k.
\end{equation}
By explicitly solving for $V_m(z)$, we obtain the following integral representation for $G(k,m)$.
\begin{theorem}\label{th1}
The number of reachable configurations $G(k,m)$ has the following exact expression
\begin{equation}\label{s2_th1_G}
G(k,m)=\frac{1}{2\pi i}\oint_\mathcal{C}(-1)^mz^{m-k-1}V_m(z)dz.
\end{equation}
Here $\mathcal{C}$ is a closed counterclockwise contour around the origin in the $z$-plane, with $|z|<1/2$ on $\mathcal{C}$, 
\begin{equation}\label{s2_th1_V}
V_m(z)=\frac{z^{1+m(m+1)/2}}{S(z)}\sum_{n=1}^\infty (-1)^{n+m}z^{n(n+1)/2+nm}\prod_{L=0}^m \frac{1}{1-z^{L+n}}\prod_{L=1}^{n-1}\frac{1}{(1-z^L)^2}
\end{equation}
and 
\begin{equation}\label{s2_th1_S}
S(z)=(2z^2-3z+2)S_1(z)-(4z^2-4z+1)S_2(z)+2z^2-z-1,
\end{equation}
where
\begin{equation}\label{s2_th1_Sk}
S_k(z)=\sum_{i=1}^\infty (-1)^{i+1}z^{i^2/2+(2k-1)i/2}\prod_{j=1}^i\frac{1}{(1-z^j)^2}.
\end{equation}
\end{theorem}

Since the total number of reachable configurations for the original problem is $G(k)=G(k,0)$, using (\ref{s2_th1_G}) in (\ref{s2_G0ic}) with $m=1$ we get the exact expression of $G(k)$:
\begin{corollary}
An exact expression for the total number of reachable configurations $G(k)$ is
\begin{equation*}\label{s2_co1_G}
G(k)=2^{k-2}+\frac{1}{2\pi i}\oint_\mathcal{C}\frac{2^k-z^{-k}}{1-2z}V_1(z)dz,
\end{equation*}
where $V_1(z)$ is given by (\ref{s2_th1_V}) with $m=1$.
\end{corollary}

Using (\ref{s2_th1_G}), we obtain the following asymptotic expressions for $G(k,m)$. These asymptotic results were first obtained in \cite{Knessl}, but there some parameters could only be determined numerically. 
\begin{theorem}\label{th2}
The number of reachable configurations $G(k,m)$ has the following asymptotic expressions:
\begin{enumerate}
\item When $k\to\infty$ and $m=O(1)$,
\begin{equation}\label{s2_th2_a}
G(k,m)\sim \frac{z_*^{m(m+3)/2-k}}{S'(z_*)}\sum_{n=1}^\infty (-1)^{n+1}z_*^{n(n+1)/2+nm}\prod_{L=0}^m \frac{1}{1-z_*^{L+n}}\prod_{L=1}^{n-1}\frac{1}{(1-z_*^L)^2},
\end{equation}
where $z_*>0$ is the unique root of $S(z)=0$ for $|z|<1/2$ and $S'(z_*)$ is the first derivative of $S(z)$ evaluated at $z=z_*$. The numerical approximation of $z_*$ to 15 decimal places is given below
$$z_*=0.43072\;95931\;37930\cdots.$$

\item When $k,\,m\to\infty$ with $2k/m^2>1$,
\begin{equation}\label{s2_th2_b}
G(k,m)\sim\frac{z_*^{m(m+5)/2-k+1}}{S'(z_*)}\prod_{L=0}^\infty\frac{1}{1-z_*^{L+1}}.
\end{equation}

\item When $l=k-m(m+5)/2$ and $2\le l\le m+3$,
\begin{equation}\label{s2_th2_c1}
G(k,m)=-\frac{1}{(l-2)!}\lim_{z\to 0}\frac{d^{l-2}}{dz^{l-2}}\bigg[\frac{1}{S(z)}\prod_{L=0}^\infty\frac{1}{1-z^{L+1}}\bigg].
\end{equation}
When $k$ and $m\to\infty$ with $l>m+3$, (\ref{s2_th2_c1}) holds asymptotically, and then reduces to (\ref{s2_th2_b}).

\end{enumerate}
\end{theorem}

Now we consider the asymptotic expression for $G(k)$ when $k\to\infty$. Since $z_*<1/2$, we can neglect the term $2^{k-2}$ in the right of (\ref{s2_G0ic}) and use (\ref{s2_th2_a}) in (\ref{s2_G0ic}) with $m=1$. Since $z_*<1/2$, we have
\begin{equation*}
\sum_{l=1}^k 2^{k-l}z_*^{-l}\sim\frac{z_*^{-k}}{1-2z_*},\; k\to\infty,
\end{equation*}
and hence the following corollary.
\begin{corollary}
The total number of reachable configurations $G(k)$ has the following asymptotic expression, when $k\to\infty$:
\begin{equation*}\label{s2_co2_G}
G(k)\sim \frac{z_*^{2}}{(1-2z_*)S'(z_*)}\bigg[\sum_{n=1}^\infty\frac{ (-1)^{n+1}z_*^{n(n+3)/2}}{(1-z_*^n)(1-z_*^{n+1})}\prod_{L=1}^{n-1}\frac{1}{(1-z_*^L)^{2}}\bigg]\bigg(\frac{1}{z_*}\bigg)^k.
\end{equation*}
Here $S(z)$ is as in (\ref{s2_th1_S}) and (\ref{s2_th1_Sk}). 
\end{corollary}

\section{Brief derivations}
Using the generating function (\ref{s2_Vm}) in (\ref{s2_Gk}), we obtain the following recurrence equation for $V_m(z)$
\begin{equation*}\label{s3_V_recu}
V_{m+1}(z)+V_{m-1}(z)=(2-z^{-m-1})V_m(z).
\end{equation*}
We notice that this equation is of the same form as (3.12) given in \cite{Knessl}, with $z=1/a$. Two linearly independent solutions are given by (3.29) and (3.30) in \cite{Knessl}. We reject the growing solution given by (3.30) in \cite{Knessl} since we expect that $V_m(z)$ will be bounded as $m\to\infty$. Hence, we have 
\begin{eqnarray}
V_m(z)&=&C(z)\frac{2\pi}{\log z}\exp\bigg(\frac{\pi^2}{2\log z}\bigg)\sum_{J=m+1}^\infty (-1)^{m+J+1}z^{-(m-J)^2/2-(m-J)/2}\nonumber\\
&&\times\prod_{L=0}^m\frac{1}{1-z^{L-J}}\prod_{L=m+1,\,L\ne J}^\infty\frac{1}{(1-z^{L-J})^2}\nonumber\\
&=&C(z)\frac{2\pi}{\log z}\exp\bigg(\frac{\pi^2}{2\log z}\bigg)\sum_{n=1}^\infty (-1)^{n+1}z^{-n^2/2+n/2}\prod_{L=0}^m\frac{1}{1-z^{-L-n}}\nonumber\\
&&\times\prod_{L=1,\,L\ne n}^\infty\frac{1}{(1-z^{L-n})^2}\label{s3_V_1},
\end{eqnarray}
where $C(z)$ is independent of $m$.
We note that we can rewrite the two products in (\ref{s3_V_1}) as follows:
\begin{equation}\label{s3_V_prod1}
\prod_{L=0}^m\frac{1}{1-z^{-L-n}}=(-1)^{m+1}z^{(n+m/2)(m+1)}\prod_{L=0}^m \frac{1}{1-z^{L+n}}
\end{equation}
and 
\begin{equation}\label{s3_V_prod2}
\prod_{L=1,\,L\ne n}^\infty\frac{1}{(1-z^{L-n})^2}=z^{n(n-1)}\prod_{L=1}^{n-1}\frac{1}{(1-z^L)^2}\prod_{L=1}^\infty\frac{1}{(1-z^L)^2}.
\end{equation}
Using (\ref{s3_V_prod1}) and (\ref{s3_V_prod2}) in (\ref{s3_V_1}) leads to
\begin{eqnarray}
V_m(z)&=&C(z)\frac{2\pi}{\log z}\exp\bigg(\frac{\pi^2}{2\log z}\bigg)\prod_{L=1}^\infty\frac{1}{(1-z^L)^2}\bigg\{\sum_{n=1}^\infty (-1)^{n+m}z^{n^2/2-n/2+(n+m/2)(m+1)}\nonumber\\
&&\times\bigg[\prod_{L=0}^m\frac{1}{1-z^{L+n}}\bigg]\bigg[\prod_{L=1}^{n-1}\frac{1}{(1-z^L)^2}\bigg]\bigg\}.\label{s3_V_2}
\end{eqnarray}
To determine $C(z)$, we use (\ref{s2_Vm}) in the boundary condition (\ref{s2_G1ic}), which yields
\begin{equation}\label{s3_V_C}
\Big(z^4+2z^2-1+\frac{z^3}{1-2z}\Big)V_1(z)-z^2V_2(z)=\frac{z^4}{1-2z}.
\end{equation}
We note that (\ref{s3_V_C}) holds for $|z|<1/2$. We introduce the function $U_m(z)$
\begin{equation*}\label{s3_Um}
U_m(z)=\sum_{n=1}^\infty (-1)^{n+m}z^{n^2/2-n/2+(n+m/2)(m+1)}\prod_{L=0}^m\frac{1}{1-z^{L+n}}\prod_{L=1}^{n-1}\frac{1}{(1-z^L)^2},
\end{equation*}
and then $V_m(z)$ in (\ref{s3_V_2}) is 
\begin{equation}\label{s3_V_3}
V_m(z)=C(z)\frac{2\pi}{\log z}\exp\bigg(\frac{\pi^2}{2\log z}\bigg)\bigg[\prod_{L=1}^\infty\frac{1}{(1-z^L)^2}\bigg]U_m(z).
\end{equation}
Using (\ref{s3_V_3}) with $m=1$ and $m=2$ in (\ref{s3_V_C}), we obtain $C(z)$ as
\begin{eqnarray}\label{s3_C_u}
C(z)&=&\Big[(-2z^5+z^4-3z^3+2z^2+2z-1)U_1(z)-z^2(1-2z)U_2(z)\Big]^{-1}\nonumber\\
&&\times\frac{\log z}{2\pi}\exp\bigg(-\frac{\pi^2}{2\log z}\bigg)z^4\prod_{L=1}^\infty(1-z^L)^2.
\end{eqnarray}
Instead of using (\ref{s3_C_u}) in (\ref{s3_V_3}), we introduce the functions $S_k(z)$ defined in (\ref{s2_th1_Sk}) and rewrite $C(z)$ in terms of the $S_k(z)$. (The reason for making this change is that it allows us to verify the equivalence of $S(z)=0$ with equation (3.39) in \cite{Knessl}, which we will discuss later.) Using (\ref{s2_th1_Sk}), $U_1(z)$ and $U_2(z)$ in (\ref{s3_C_u}) can be expressed as 
\begin{equation}\label{s3_U1}
U_1(z)=-S_1(z)+\Big(1+\frac{1}{z}\Big)S_2(z)-\frac{1}{z}S_3(z)
\end{equation}
and
\begin{equation}\label{s3_U2}
U_2(z)=-S_1(z)+\Big(1+\frac{1}{z}+\frac{1}{z^2}\Big)S_2(z)-\Big(\frac{1}{z}+\frac{1}{z^2}+\frac{1}{z^3}\Big)S_3(z)+\frac{1}{z^3}S_4(z).
\end{equation}
We rewrite $S_k(z)$ in (\ref{s2_th1_Sk}) as
\begin{equation*}\label{}
S_k(z)=\sum_{i=1}^\infty (-1)^{i+1}z^{i^2/2+(2k-1)i/2}\prod_{j=1}^{i+1}\frac{1}{(1-z^j)^2}(1-2z^i+z^{2i}),
\end{equation*}
and then, after some calculation, we obtain the following recurrence equation
\begin{equation}\label{s3_Sk_recu}
z^{1-k}S_{k-1}(z)+(1-2z^{1-k})S_k(z)+z^{1-k}S_{k+1}(z)=1.
\end{equation}
Hence, we can use (\ref{s3_Sk_recu}) with $k=2$ and $k=3$ to eliminate $S_3(z)$ and $S_4(z)$ from (\ref{s3_U1}) and (\ref{s3_U2}). Thus, after some simplification, we obtain $C(z)$ as
\begin{equation}\label{s3_C_s}
C(z)=\frac{z}{S(z)}\frac{\log z}{2\pi}\exp\bigg(-\frac{\pi^2}{2\log z}\bigg)\prod_{L=1}^\infty(1-z^L)^2,
\end{equation}
where $S(z)$ is defined in (\ref{s2_th1_S}).
Using (\ref{s3_C_s}) in (\ref{s3_V_3}), we obtain (\ref{s2_th1_V}) in Theorem 2.1. Then we have the integral expression of $G(k,m)$ in (\ref{s2_th1_G}). 

In the remainder of this section, we discuss the asymptotic approximations to $G(k,m)$. We first consider $k\to\infty$ and $m=O(1)$. By using Rouch\'e's theorem, or by plotting $S(z)$ given in (\ref{s2_th1_S}) numerically, we notice that there is a single real root for $|z|<1/2$. We denote this single root as $z_*$, which satisfies $S(z_*)=0$. Then in the $z$-plane with $|z|<1/2$, the integrand in (\ref{s2_th1_G}) has a simple pole at $z=z_*$, which is the dominant singularity. We use the residue theorem to evaluate the integral in (\ref{s2_th1_G}) asymptotically, which yields (\ref{s2_th2_a}). We note that it's not hard to verify that $S(z)=0$ is equivalent to (3.39) in \cite{Knessl}, after we substitute $1/z$ for $a$ in (3.39). Thus $z_*^{-1}=a$, whose numerical approximation to 100 decimal places is given in \cite{Knessl}. By comparing (\ref{s2_th2_a}) in this paper with formula (3.37) in \cite{Knessl}, we obtain an analytic expression for the constant $c_1$ in (3.38) in \cite{Knessl} as follows:
\begin{equation}\label{s3_c1}
c_1=-\frac{\log z_*}{2\pi}\exp\bigg(-\frac{\pi^2}{2\log z_*}\bigg)\frac{1}{S'(z_*)}\prod_{j=1}^\infty(1-z_*^j)^2.
\end{equation}
By evaluating (\ref{s3_c1}) numerically, we find that the numerical values of $c_1$ and $c_1K_*$ provided in (3.42) and (3.43) in \cite{Knessl} are only correct to about 5 decimal places, though they are given to 15 places. The correct values to 15 decimal places are
$$c_1=2.02740\; 20474\; 68498\cdots$$
and 
$$c_1K_*=0.28777\; 77049\; 35052\cdots.$$

In the limit when $k,\,m\to\infty$ simultaneously, the $n=1$ term in (\ref{s2_th2_a}) dominates, which yields (\ref{s2_th2_b}). This recovers the result given in (4.18) in \cite{Knessl}, but now with $c_1$ computed explicitly.

Next, we consider $k$ and $m$ large with $l=k-m(m+5)/2=O(1)$. Note that $l=2$ corresponds to the smallest number of possible configurations, as $G(k,m)=0$ for $l\le 1$. Following \cite{Knessl} we denote $G(k,m)$ as $W(l,m)$ and rewrite (\ref{s2_th1_G}) as 
\begin{eqnarray}
G(k,m)\equiv W(l,m)&=&\frac{1}{2\pi i}\oint_\mathcal{C}\frac{1}{z^{l-1}}\frac{1}{S(z)}\sum_{n=1}^\infty (-1)^nz^{(n-1)(n/2+m+1)}\nonumber\\
&&\times\prod_{L=0}^m\frac{1}{1-z^{L+n}}\prod_{L=1}^{n-1}\frac{1}{(1-z^L)^2}dz.\label{s3_G}
\end{eqnarray}
For a sufficiently small closed contour $\mathcal{C}$ around the origin in the $z$-plane and $l\ge 2$, $z=0$ is the only pole inside of $\mathcal{C}$, and it is of order $l-1$. Thus, using the residue theorem in (\ref{s3_G}) leads to 
\begin{eqnarray}
W(l,m)&=&\frac{1}{(l-2)!}\lim_{z\to 0}\frac{d^{l-2}}{dz^{l-2}}\bigg\{\frac{1}{S(z)}\sum_{n=1}^\infty (-1)^nz^{(n-1)(n/2+m+1)}\nonumber\\
&&\times\prod_{L=0}^m\frac{1}{1-z^{L+n}}\prod_{L=1}^{n-1}\frac{1}{(1-z^L)^2}\bigg\}.\label{s3_th2_b}
\end{eqnarray}
The first two terms from the infinite sum in the right-hand side of (\ref{s3_th2_b}) are
\begin{equation}\label{s3_th2_b2term}
-\prod_{L=0}^m\frac{1}{1-z^{L+1}}+\frac{z^{m+2}}{(1-z)^2}\prod_{L=0}^m\frac{1}{1-z^{L+2}},
\end{equation}
and the remaining terms are of order $O(z^{2m+5})$ as $z\to 0$. 
We notice that after taking $l-2$ derivatives, the second term in (\ref{s3_th2_b2term}) is of order $O(z^{m-l+4})$ as $z\to 0$. This implies that as long as $m-l+4\ge 1$, i.e., $l\le m+3$, only the $n=1$ term in the infinite sum in (\ref{s3_th2_b}) contributes to the derivative at $z=0$, which yields
\begin{equation}\label{s3_th2_b2}
W(l,m)=-\frac{1}{(l-2)!}\lim_{z\to 0}\frac{d^{l-2}}{dz^{l-2}}\bigg[\frac{1}{S(z)}\prod_{L=0}^m\frac{1}{1-z^{L+1}}\bigg],\; 2\le l\le m+3.
\end{equation}
Rewriting the expression in the brackets in (\ref{s3_th2_b2}) as
\begin{equation*}\label{s3_th2_b2prod}
\frac{1}{S(z)}\prod_{L=0}^m\frac{1}{1-z^{L+1}}=\frac{1}{S(z)}\prod_{L=0}^\infty\frac{1}{1-z^{L+1}}\prod_{L=m+1}^\infty(1-z^{L+1}),
\end{equation*}
we see that (\ref{s3_th2_b2}) is independent of $m$. This leads to (\ref{s2_th2_c1}), and gives an analytic expression for $W_0(l)$, which appears in \cite{Knessl} but there no explicit expression is given. This concludes our derivation.

\end{document}